\documentclass[a4paper,12pt]{article}
\usepackage{booktabs}
\usepackage{float}
\usepackage{lineno}
\modulolinenumbers[5]
\usepackage[intlimits]{amsmath}
\usepackage[algoruled,resetcount]{algorithm2e}
\usepackage{amssymb,amsthm,amsfonts}
\usepackage{graphicx}
\usepackage{setspace}
\usepackage[center]{subfigure}
\usepackage{paralist}
\usepackage{color}
\definecolor{cor}{rgb}{0,0,0}
\definecolor{cor1}{rgb}{0,0,0}
\usepackage[abs]{overpic}

\newcommand{\mesh}{{\{$\mathbb T_l\}_{0 \leq l \leq L}$}}
\newcommand{\ti}{TOL_{iter}}
\newcommand{\tn}{TOL_{newton}}
\newcommand{\ri}{res_{iter}}
\newcommand{\rn}{res_{newton}}
\newcommand{\maxit}{{MAX}_{iter}}
\newcommand{\refAB}{\eqref{eq:AB}}

\newcommand{\VV}{V^o}
\newcommand{\VU}{V^p}
\newcommand{\VUH}{V^p_h}

\newcommand{\sums}{\Sigma s}
\newcommand{\sumn}{\Sigma n}

\newcommand{\V}{V}

\newcommand{\BB}{{C}}

\newcommand{\NO}{{N^o}}
\newcommand{\NP}{{N^p}}

\newcommand{\A}{{\mathcal{A}}}
\newcommand{\B}{{\mathcal{B}}}

\newcommand{\R}{{\mathbb{R}}}

\newcommand{\abs}[1]{ \left\vert #1 \right\vert} 

\newtheorem{remark}{Remark}[section]

\newcommand{\newln}{\\&\quad\quad{}}

\newcommand{\grd}{{\{$\mathbb T_l\}_{0 \leq l \leq L}$}}

\newcommand{\vnoti}{{w^0_i}}
\newcommand{\vonei}{{w^1_i}}

\newcommand{\kd}{{k_d}}
\newcommand{\kon}{{k_{on}}}
\newcommand{\kono}{{k_{on}}}
\newcommand{\krec}{{k_{rec}}}
\newcommand{\kdeg}{{k_{deg}}}
\newcommand{\kir}{{k_{iR}}}
\newcommand{\kib}{{k_{iB}}}
\newcommand{\kic}{{k_{iC}}}
\newcommand{\koff}{{k_{off}}}
\newcommand{\koffo}{{k_{off}}}
\newcommand{\kk}{{K}}
\newcommand{\dd}{{\mu}}
\newcommand{\ut}{{\tilde{u}}}

\newcommand{\du}{u}
\newcommand{\dv}{v}

\usepackage[unicode, bookmarks, colorlinks, breaklinks,
	pdftitle={Carraro},
	pdfauthor={Carraro},
	pdfproducer={}
]{hyperref}  
\hypersetup{linkcolor=blue,citecolor=blue,filecolor=black,urlcolor=blue}

\begin{document}
\title{Coupling vs decoupling approaches for PDE/ODE systems modeling intercellular signaling}

\author{Thomas Carraro\thanks{thomas.carraro@iwr.uni-heidelberg.de} \and ~Elfriede Friedmann \and ~Daniel Gerecht\vspace{0.5cm}\\
  Institute for Applied Mathematics\\ Heidelberg University, Germany
}

\maketitle

\begin{abstract}
We consider PDE/ODE systems for the simulation of intercellular signaling in multicellular environments. 
{\color{cor}The intracellular processes for each cell described here by ODEs determine the long-time dynamics, but the PDE part dominates the solving effort.}
Thus, it is not clear if commonly used decoupling methods can outperform a coupling approach. Based on a sensitivity analysis, we present a systematic comparison between coupling and decoupling approaches for this class of problems and show numerical results.
{\color{cor}For biologically relevant configurations of the model, our quantitative study shows that a coupling approach performs much better than a decoupling one}.

\end{abstract}
\emph{keywords}
Coupled PDE/ODE systems, Sensitivity analysis, Multilevel preconditioner, Intercellular Signaling

\section{Introduction}
\label{intro}
{\color{cor}
Cellular signaling has been mathematically described by a variety of models mostly relying on large systems of ordinary differential equations (ODE) \cite{kestler2008network}. These earlier models were extended by partial differential equations (PDE) to accurately consider concentration gradients and their effect \cite{sherratt1995receptor,marciniak2003receptor,friedmann2013interaction,Busse:2010}. In our intercellular model we consider the diffusion and action of small signaling  proteins in the intercellular space described by PDE, i.e. reaction-diffusion equations, coupled on the cell surfaces with ODEs responsible for the intracellular dynamics.
This coupling of mixed differential equations allows mainly two strategies for an implicit solver:
}
(1) nonlinear methods among them the nonlinear multigrid method  also called ``full approximation scheme'' (FAS) \cite{Brandt:1977, Hackbusch:1985},
(2) linearization based approaches (Newton-type).
These methods can be used in a combined approach, where for example a Newton-type method can be used as smoother for a FAS and a linear or a nonlinear multigrid can be used as a preconditioner for a Newton-type method.
The comparison and discussion of advantages and disadvantages of these strategies that depend on many aspects like, e.g.\ the accuracy of the Jacobian approximation \cite{Mavriplis:2002}, is not the focus of our work.
Moreover, since in the considered coupled PDE/ODE system the linearization is not a critical point we  choose a Newton-type method preconditioned by a linear multigrid and study the effect of splitting the linearization. A decoupling solution approach, {\color{cor}based on a fixed-point method}, is often used when restrictions on accuracy can be relaxed in order to allow an easier numerical treatment of complicated problems.
Such an approach makes it possible to reuse existing solvers and is widely used in numerical methods for coupled systems, see \cite{EsmailyMoghadam201363,farhat1998fast,felippa2001partitioned,mok2001partitioned,heil1998stokes,heuveline_ode}.
{\color{cor}
In case of strongly coupled equations (see Section \ref{section:coupling} for the definition of \emph{strong} coupling used here), this strategy leads usually to high computational costs through very small time steps needed to reduce the strength of the coupling.
In fact, the convergence of the fixed-point iterations is typically linear with the convergence rate depending on the used block-iterative method and on the strength of the coupling since it depends on the spectral radius of the matrices involved \cite{Cervera:1996}.
Therefore a drawback of decoupling solvers is their low convergence and possibly divergence.
On the contrary, the drawback of fully implicit solvers is mainly that they demand for the implementation of a special-purpose code.
}

We consider PDE/ODE systems for the simulation of intercellular signaling in multicellular environments. Since the ODE part does not lead to a large discretization system like the PDE part, it is not clear if a decoupling method can outperform a coupling approach. {\color{cor}In fact other works with PDE/ODE models have used decoupling approaches for systems arising in biological applications, e.g.\ blood flow including chemical interaction \cite{Ricken:2015, Quarteroni:2003}, signal transduction \cite{friedmann2013interaction}, cardiovascular flow \cite{EsmailyMoghadam201363}, cancer invasion \cite{Surulescu:2014}.
{\color{cor1} With the exception of \cite{EsmailyMoghadam201363} these works do not consider a comparison with a monolithic approach}. In addition, for {\color{cor1} the most of} these applications it is not clear whether the strength of the coupling is strong or weak.}

In this context, the scope of our work is to present a systematic comparison between coupling and decoupling approaches for this class of problems.
The method is based on a sensitivity analysis to compute the strength of the coupling.
Additionally, we compare a multigrid method in which the coupling is considered only at the coarsest level to a fully coupling approach.
There are few works that deal with the fully coupling solution process of dimensionally heterogeneous systems such as \cite{Leiva:2010}.
Furthermore, we focus on the solution of local microenvironments. Therefore, this solution process can be used for example as local solver for nonlinear preconditioner of Newton-type methods \cite{Cai:2002} or domain decomposition methods \cite{quarteroni1999domain}.

\paragraph{Outline}
The paper is organized as follows. In Section \ref{sec:problem} we give an abstract description of the model. We present the mathematical formulation and the functional setting. We discretize the coupled PDE/ODE system by the finite element method (FEM) in Section \ref{section:discretization}. We use a sensitivity approach in Section \ref{section:coupling} to analyze the coupling of the PDE/ODE system and in Section \ref{sec:numerical schemes} we present different solving approaches for the coupled system. We present the numerical results exemplarily for a particular application and discuss numerical aspects in Section \ref{sec:results}. In Section \ref{application} we describe a realistic configuration and give a biological interpretation to the results obtained.
\section{Mathematical Models for Intercellular Signaling}\label{sec:problem}
\begin{figure}[H]
  \caption{Visualization of the computational domain}
\centering
\subfigure[8 interacting cells with surfaces $\Gamma_i$]{
  \begin{overpic}[scale=0.5,unit=1mm]{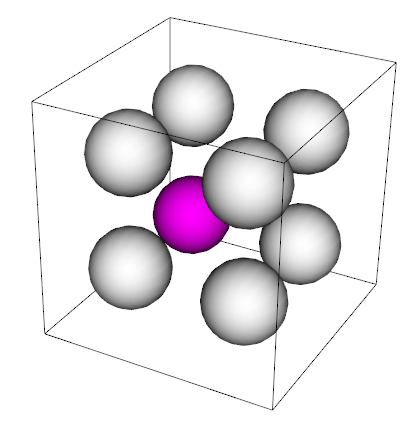}
    \put(-3,28){$\Gamma_{out}$}
    \put(32,45){$\Gamma_{i}$}
  \end{overpic}
}
\subfigure[Intercellular area $\Omega$ \newline(sliced for visualization)]{
  \includegraphics[width=.47\textwidth]{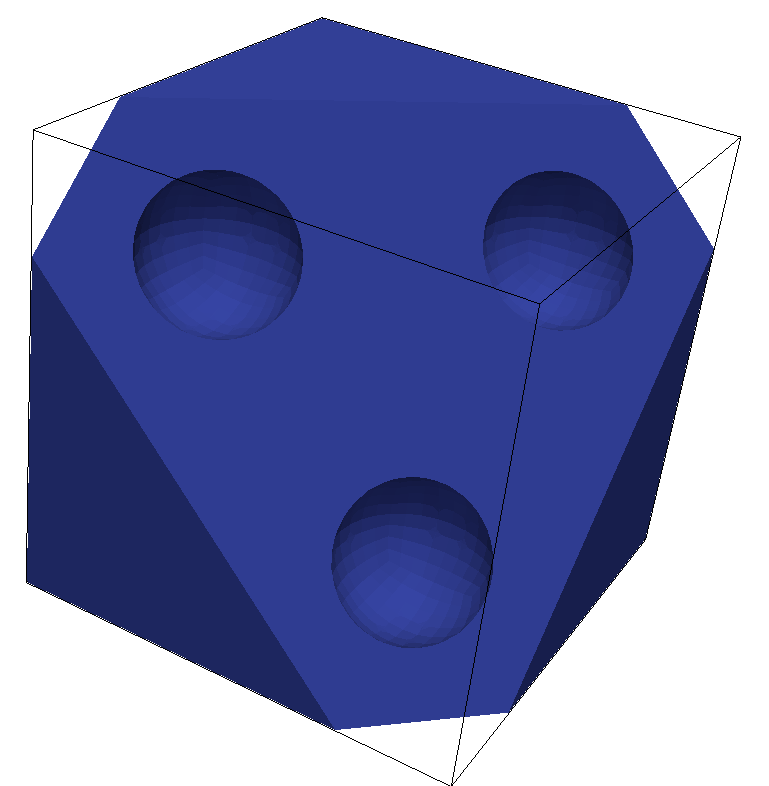}}

\label{3d_pic2}
\end{figure}

{\color{cor}
Our intercellular signaling model consists of one PDE equation for the interaction between the cells in the intercellular area $\Omega \subset \R^3$ coupled with ODEs for the intracellular processes. 
}
We denote by $N_c$ the number of cells in $\Omega$ and indicate by $\Gamma_i$ the boundary of each cell $i$ for $i=1,\dots,N_c$. The outer boundary of $\Omega$ is denoted by $\Gamma_{out}$, see Figure \ref{3d_pic2}.

Depending on the type of intercellular signaling, different nonlinear operators describe the dynamics in the intercellular area ($\A_\Omega$), e.g.\ degradation, the dynamics on the cell surfaces ($\A_{\Gamma_i}$) of each cell and the intracellular processes ($\B_i$). We denote the solution of the PDE part with $u$ and the vector of solutions of the ODE part with $v$. 
\begin{align}
\begin{alignedat}{10}
\partial_t u(t,x) - \dd \Delta u(t,x)+ \A_{\Omega}(u(t,x))&=0\quad&& \text{for }(t,x)\in (0, T]\times\Omega, \\
\dd \partial_n u(t,x)-\A_{\Gamma_i}(u(t,x),v_i(t))&=0&& \text{for }(t,x) \in (0,T]\times\Gamma_i , \\
\dd \partial_n u(t,x)&=0&& \text{for }(t,x) \in (0, T]\times\Gamma_{out} , \\
\partial_t v_i(t) + \B_i(\ut_i(t),v_i(t))&=0\label{eq:coupling} &&\text{for } t \in(0, T],\\
\end{alignedat}
\end{align}
with given initial values $u(0,x)=u^0 $ and $v(0)=v^0$. We denote the average of $u$ on the surface of $\Gamma_i$ by $\ut_i$ and by $v_i$ the associated ODE values with this cell
\begin{align}
  \label{u tilde}
\ut_i(t)&=\frac{\int_{\Gamma_i}u(t,s)\text{ ds}}{\abs{\Gamma_i}}.
\end{align}
{\color{cor}
}

\begin{remark}
\label{remark:stationary}
To study the dynamical process and validate the model we compute the entire trajectory. 
{\color{cor}
Nevertheless, the simulations converge to a stable steady state. Therefore, we consider as well a coupling and a decoupling solver for a computation of the steady state in Section \ref{section:stat_results}.
}
\end{remark}
\section{Discretization}
\label{section:discretization}
For a variational formulation we introduce the Hilbert space $\VU=H^1$, for the PDE part of the equation, and the vector space $\VV=\R^{n}$, where $n$ denotes the number of ordinary differential equations in the system. We define the product space $\V:=\VU\times \VV$.

We consider  the implicit Euler method as time stepping scheme, and spatially discretize the computational domain $\Omega$ by continuous finite elements.

Considering a time step $k$ we use the semi-discretized weak formulation of the problem \eqref {bio_equation} to compute $(u^{n+1},v^{n+1})\in \V$ in each time step for all $\varphi\in \VU$:
\begin{align}\label{eq:semiproblem}
\begin{split}
(u^{n+1}, \varphi)_\Omega + k\mu(\nabla u^{n+1},\varphi)+ k (\A_\Omega(u^{n+1}), \varphi)_\Omega\\+ k\sum_{i\leq N_C} (\A_{\Gamma_i}(u^{n+1},v^{n+1}_i), \varphi)_{\Gamma_i} &=  (u^n,\varphi)_\Omega,
\\
 v^{n+1} + k \sum_{i\leq N_C}\B_i(\ut_i^{n+1},v_i^{n+1})& = v^n.
\end{split}
\end{align}

{\color{cor}The domain in which the PDE part is defined is the volume between the T cells. Figure \ref{3d_pic2} (b) depicts a portion of it showing that the T cells are holes in the domain. To discretize the system \eqref{eq:semiproblem}, we define a grid \mesh, consisting of non-overlapping hexahedral cells $K$. The volume of the holes, representing the T cells, is not discretized because the T cells are described by the ODE part of the system. Therefore the domain has an outer boundary $\Gamma_{out}$ and as many inner boundaries $\Gamma_i$ as the number of T cells, see Figure \ref{3d_pic2} (a).
The diameters of the hexahedral cells $h_K$ define a mesh parameter $h$ by the piece-wise constant function $h_{|K} = h_K$. We use the index $h$ to indicate all terms of the equations that are discretized using the grid $\mathbb T$.
}

The discrete solution component $u_h$ is sought in the finite dimensional space $\VUH \subset \VU$. We choose $\VUH$ as the space of  $Q^1$-elements, the space of functions obtained by transformations of {\color{cor}trilinear} polynomials defined on a reference unit cell {\color{cor} (regular hexahedral)}, see e.g.\ \cite{ciarlet:1978} for an introduction to the finite element method.
{\color{cor}As explained above, only the PDE part needs to be discretized by the FEM, but due to {\color{cor}the coupling terms $\tilde u_i$, see \eqref{u tilde},} the ODE part of the discretized system depends as well on the spatial discretization. Therefore, we use the symbol $v_h$ to indicate this dependency. Since $v$ was already defined on a finite dimensional space $\VV$, we have $v_h \in \VV$.}
Then we can write the fully discretized version of system \eqref{eq:semiproblem} as follows for all $\varphi_h\in \VUH$:
\begin{align}\label{eq:instationary_problem}
\begin{split}
(u_h^{n+1}, \varphi_h)_\Omega + k (\A_\Omega(u_h^{n+1}), \varphi_h)_\Omega+k \mu (\nabla u_h^{n+1},\nabla\varphi_h)\\+ k\sum_{i\leq N_C} (\A_{\Gamma_i}(\ut_{h,i}^{n+1},v_{h,i}^{n+1}), \varphi_h)_{\Gamma_i} &=  (u^n_h,\varphi_h)_\Omega, 
\\
v_h^{n+1} + k \sum_{i\leq N_C}\B_i(\ut_{h,i}^{n+1},v_{h,i}^{n+1}) &= v^n_h.
 \end{split}\end{align}

A nonlinear coupled system arises as  discrete system \eqref{eq:instationary_problem} and needs to be solved for each step of the time marching scheme. For the resulting discrete nonlinear coupled system we introduce the shorter notation 
\begin{align}
\tag{$4'$}\label{eq:AB}
\begin{split}
 A_h (u,v)=f ,
\\
 B_h (u,v)=g .
\end{split}
\end{align}
We use the subscript $h$ to indicate the dependence of the operator $B_h$ on the mesh discretization through the coupling with the PDE part. We omit the subscript $h$ for the solution components $u$ and $v$ to {\color{cor1}simplify} the notation in the next sections.

\section{Sensitivity Analysis of the Coupled System}
\label{section:coupling}
{\color{cor} 
In this section we study the strength of the coupling between the PDE and the ODE part of the system.
This is done to motivate the choice of the solution process that we present in the following part of the work.
In this work we refer to \emph{strong} or \emph{weak} coupling meaning the \emph{strength} of the coupling independently of the approach used to solve the coupled problem. We first define in the next subsection the used definition of \emph{strong} and \emph{weak} coupling, then we apply this definition to our system of equations. 
}

\subsection{A criterion to measure the coupling strength}

Let $S_h: v \mapsto u$ and $T_h: u \mapsto v$ denote the solution operator for the PDE part and respectively for the ODE part of the discretized system of equations \refAB.
The first equation, $u = S_h(\dv)$, is solved for a given value of $v$, then the second equation, $v = T_h(u)$, is solved with the resulting value of $u$ and the cycle is iterated until a given tolerance is reached.
This process can also be written as a composition of the two operators:
\begin{align}
\label{eq:fixed point composition}
u^{n+1} &= S_h\bigl(T_h(u^n)\bigr).
\end{align}

A fixed-point iteration to solve a coupled system of equations has a slow convergence rate (typically only linear) and the number of fixed-point iterations depends on the nature of the coupling and the model parameters.
Considering the formulation \eqref{eq:fixed point composition} to solve our coupled system we write the Jacobian of the fixed-point operator as
\begin{align}
\label{Jacobian}
J = \frac{\partial S_h}{\partial \dv} \frac{\partial T_h}{\partial \du}.
\end{align}
For convergence the fixed-point iteration \eqref{eq:fixed point composition} has to fulfill the following criterion according to the Banach fixed-point theorem
\begin{align}
\|J\| < 1,
\end{align}
in some norm $\|\cdot\|$. A more convenient criterion is the substitution of the norm with the spectral radius of   the matrix $J$
\begin{align}
  |\lambda_{max} (J)| < 1.
\end{align}

Based on this statement Haftka and coworkers have defined in \cite{Haftka} a quantitative measure for the strength of the coupling between the two parts of the problem.
In fact, according to their definition a system is \emph{weakly coupled} if
\begin{align}
  \label{strength criterion weak}
  |\lambda_{max}(J)| \ll 1, 
\end{align}
respectively is \emph{strongly coupled} if
\begin{align}
  \label{strength criterion strong}
 |\lambda_{max}(J)| \gg 1.
\end{align}
In \cite{Haftka} there is no precise separation of the two ranges of \emph{weak} and \emph{strong} coupling. We also do not give a quantitative definition of the ranges, because a more general definition should consider many models.
The criterion \eqref{strength criterion weak}, respectively \eqref{strength criterion strong}, has been used here with the maximal eigenvalue of \eqref{Jacobian} to define whether the coupling of our system is \emph{weak}, respectively \emph{strong}, considering that in one case the maximal eigenvalue is two orders of magnitude smaller than 1 and in the other case almost one order of magnitude larger than 1, see results in Section \ref{sec:results}.

{\color{cor1} We remark that in case of nonlinear problems the strength of the coupling could be time dependent. In particular, this is the case if the coefficients of the operator are time dependent. In our case, the coefficients are constant.
}

\subsection{Application to the PDE/ODE coupled system of equations}
In the rest of this section we present the sensitivity approach needed to calculate the largest eigenvalue of the Jacobian $J$ of the fixed-point problem, see \eqref{Jacobian}.
Therefore, we differentiate the discretized operators $A_h$ and $B_h$ and obtain the sensitivity equations
\begin{align}
\label{sens_u}
A_{h,u}^\prime(\hat \du, \hat \dv)u_{ \delta v} + A_{h,v}^\prime(\hat \du, \hat \dv)\delta v=0, \quad \forall \delta v \in \VV,
\end{align}
and
\begin{align}
\label{sens_v}
B_{h,v}^\prime(\hat \du, \hat \dv)v_{ \delta u} + B_{h,u}^\prime(\hat \du, \hat \dv)\delta \du=0, \quad \forall \delta \du \in \VUH.
\end{align}
where we have used the notation
\begin{align*}
u_{\delta v}:=\frac{\partial u }{ \partial v}( \delta v),\quad\quad\quad v_{\delta u}=\frac{\partial v }{ \partial u}( \delta u)
\end{align*}
for the sensitivities. In the decoupled system, $u_{\delta v}$ indicates the variation of the PDE solution perturbing the solution of the ODE system and equivalently $v_{\delta u}$ is the variation of the ODE system for a perturbation of the PDE system.
For nonlinear systems of equations the sensitivity analysis depends on a given point of linearization $({\hat{u}},{\hat{v}})$. We compute an approximate numerical solution of the system \refAB { for} characteristic values of the parameters and choose the computed solution as point of linearization.

 Since the sensitivities in the linear solver strongly depend on the used time stepping scheme, we consider only the sensitivities for a  computation of the steady state. Then, the equations \eqref{sens_u} are stationary PDEs to be solved for each component of $\delta \dv$, while the ODE part \eqref{sens_v} consists of algebraic equations solved for each $\delta \du$. Therefore, we compute the sensitivity matrices ${\partial S_h}/{\partial \dv}$ as a $\NO \times \NP$ matrix and ${\partial T_h}/{\partial \du}$ as a $\NP \times \NO$ matrix, where $\NO$ denotes the number of ODE equations and $\NP$ the dimension of the PDE discretization.

 {\color{cor1}
The sensitivity matrices are
\begin{align}
  \frac{\partial S_h}{\partial \dv} = \setlength{\delimitershortfall}{0pt} \begin{pmatrix}
    \dfrac{\partial u }{\partial v_1}&
    \hdots&
    \dfrac{\partial u }{\partial v_i}&
    \hdots&
    \dfrac{\partial u }{\partial v_{N^o}}&
 \end{pmatrix}^T,
\end{align}
where each row $\biggr(\dfrac{\partial u }{\partial v_i}\biggl)^T$ is of the dimension of the PDE discretization and
\begin{align}
\frac{\partial T_h}{\partial \du} = \setlength{\delimitershortfall}{0pt} \begin{pmatrix}
  \dfrac{\partial v_1 }{\partial u_1} &
  \hdots&
    \dfrac{\partial v_i }{\partial u_1}&
    \hdots&
  \dfrac{\partial v_{N^o}}{\partial u_1}\\[2ex]
  \hdots & \hdots & \hdots & \hdots & \hdots\\[2ex]
  \dfrac{\partial v_1 }{\partial u_i} &
  \hdots&
    \dfrac{\partial v_i }{\partial u_i}&
    \hdots&
  \dfrac{\partial v_{N^o}}{\partial u_i}\\[2ex]
  \hdots & \hdots & \hdots & \hdots & \hdots\\[2ex]
  \dfrac{\partial v_1 }{\partial u_{N^p}} &
  \hdots&
  \dfrac{\partial v_i }{\partial u_{N^p}}&
    \hdots&
  \dfrac{\partial v_{N^o}}{\partial u_{N^p}}
 \end{pmatrix},
\end{align}
where $\dfrac{\partial v_i }{\partial u_j}$ denotes the derivatives of the $i^{\rm th}$ component of $v$ with respect to the $j^{\rm th}$ degree of freedom of $u$.
}
%
 \paragraph{{\color{cor1}Practical realization}}
In the PDE/ODE system presented in Section \ref{sec:results} the coupling between the two parts appears only at the boundaries $\Gamma_i$ and only with the first two components of $v$. Thus the product ${\partial S_h}/{\partial \dv} ~ {\partial T_h}/{\partial \du}$ decouples into a block diagonal matrix consisting of $2\times 2$ matrices for each biological cell. In addition, we need to calculate the sensitivities \eqref{sens_v} only for the restriction of $\delta \du$ on the boundaries $\Gamma_i$, which are nonetheless algebraic equations, so that the major costs to calculate the sensitivities are given by the PDE part \eqref{sens_u}.

%
\section{Numerical Schemes}
\label{sec:numerical schemes}
\begin{figure}[H]
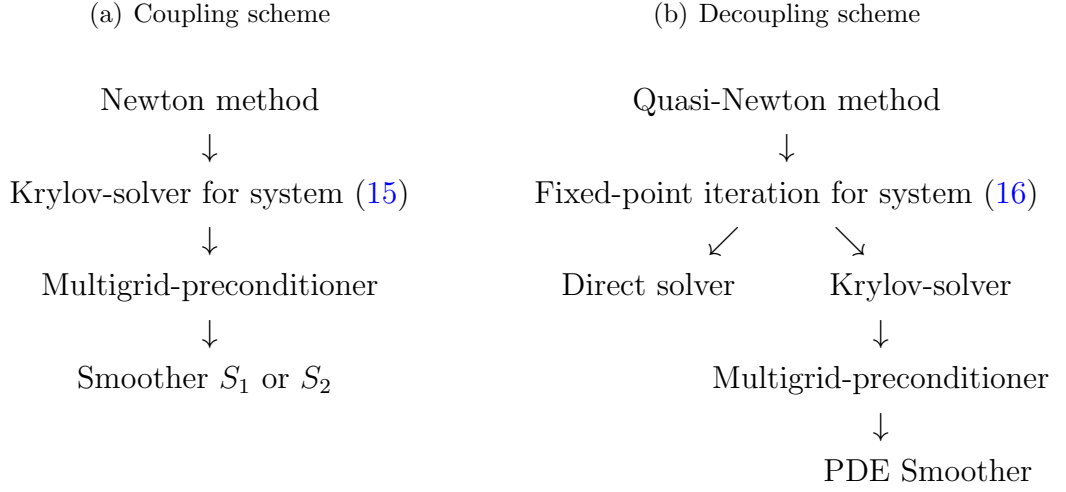

\caption{Schematic representation of the considered coupling and decoupling schemes}
\label{fig:coupling}\subfiguretopcaptrue
\subfigure[Coupling scheme]{\parbox{.5\textwidth}{\begin{gather*}
\text{Newton method}\\
\downarrow\\
\text{Krylov-solver for system \eqref{linearized problem_stat}}\\
\downarrow\\
\text{Multigrid-preconditioner}\\
\downarrow\\
\text{Smoother $S_1$ or $S_2$ }\\
\end{gather*}}}
\subfigure[Decoupling scheme]{\parbox{.5\textwidth}{\begin{gather*}
\text{Quasi-Newton method}\\
\downarrow\\
\text{Fixed-point iteration for system \eqref{linearized decoupled problem_stat}}\\
	\swarrow \quad\quad\quad		\searrow\\
\text{Direct solver}\quad\quad\quad\text{Krylov-solver}\\
\quad\quad\quad\quad\quad\quad\downarrow\\
\quad\quad\quad\quad\quad\quad\text{Multigrid-preconditioner}\\
\quad\quad\quad\quad\quad\quad\downarrow\\
\quad\quad\quad\quad\quad\quad\quad\quad\quad\text{PDE Smoother}
\end{gather*}}}
\end{figure}

In this section we present two different approaches to solve the class of coupled PDE/ODE systems presented in this work that are depicted in Figure \ref{fig:coupling}. 
In both schemes we consider a Newton-type solver for the nonlinearities. In the coupling scheme (a) the linear system is solved by a Krylov-solver {\color{cor1} that is} preconditioned by a multigrid scheme. In the decoupling scheme (b) the Newton update is approximated by a fixed-point iteration and a multigrid method is applied only to the PDE part.
Therefore, the decoupling approach is applied not to the original system but to its linearization within a Newton step.

This section is organized as follows: we first introduce the nonlinear solver based on a (Quasi-)Newton method in Subsection \ref{section:nonlinear schemes} and then we present the variants of multigrid preconditioner in Subsections \ref{mg for coupled} and \ref{mg for decoupled}.
\subsection{Nonlinear solver}
\label{section:nonlinear schemes}
 Newton-type methods provide a flexible and reliable framework for nonlinear problems by solving a series of linear equations. As explained above, we present a fully coupling and a decoupling approach to solve the linearized subproblems.

\subsubsection{Fully coupling Newton's method}
To apply Newton's method we linearize the system and solve in each Newton step the system:
\begin{align}
\label{linearized problem_stat}
 \begin{pmatrix}A_{h,u}^\prime(\du^n, \dv^n) &	 A_{h,v}^\prime(\du^n, \dv^n)\\
 B_{h,u}^\prime(\du^n, \dv^n) &  B_{h,v}^\prime(\du^n, \dv^n)\end{pmatrix}
\begin{pmatrix}
\delta \du^{n+1} \\ \delta \dv^{n+1} \end{pmatrix}
=
\begin{pmatrix}
 f- A_h(\du^n,\dv^n) \\ g - B_h(\du^n,\dv^n)
\end{pmatrix}
,
\end{align}
to obtain the Newton updates $\delta \du^{n+1}$ and $\delta \dv^{n+1}$, with which we calculate the next iterates $\du^{n+1}=\du^{n}+\delta \dv^{n}$ and $\dv^{n+1}=\dv^{n}+\delta \dv^{n+1}$. We write $A_{h,u}^\prime$ and $A_{h,v}^\prime$ for the derivatives of $A_h$ with respect to $u$ and $v$ and analogously $B_{h,u}^\prime$ and $B_{h,v}^\prime$ for the derivatives of $B_h$.

\subsubsection{Decoupling inexact Newton's method}
Secondly we consider a decoupling solving scheme for the linear systems defined in each Newton-step.
 For each Newton step $n$ the decoupled system is solved by the following fixed-point iteration
\begin{align}
\label{linearized decoupled problem_stat}
\begin{pmatrix}
 A_{h,u}^\prime(\du^n, \dv^n) &	 A_{h,v}^\prime(\du^n, \dv^n)\\
 0 &  B_{h,v}^\prime(\du^n, \dv^n)
\end{pmatrix}
\begin{pmatrix}\delta \du^{i+1} \\ \delta \dv^{i+1} 
\end{pmatrix}
=
\begin{pmatrix} f- A_h(\du^n,\dv^n) \\ g - B_h(\du^n,\dv^n) - B_{h,v}^\prime(\du^n, \dv^n) \delta  \du^{i}
\end{pmatrix}
\end{align}
 until the Newton updates  $(\delta \du^{i+1},\delta \dv^{i+1})$ fulfill the linear residual of the system \eqref{linearized problem_stat} to an accuracy ($\ti$). 

   If the linear system is solved until the full accuracy is reached then the update formula \eqref{linearized decoupled problem_stat} defines a full Newton-step, which is equivalent to the exact solution of \eqref{linearized problem_stat}.
   A common approach to accelerate the computation of the Newton-updates is a Quasi-Newton iteration in which the Jacobian matrix is approximated only up to a certain accuracy. This corresponds not to solving the system \eqref{linearized decoupled problem_stat} to the full accuracy.
   In this way the costs per Newton iteration are reduced, while the number of Newton iterations increases. A trade-off between accuracy and total costs  can enable a reduction of computing time with respect to a full Newton method. Such a Quasi-Newton scheme is obtained if a low accuracy ($\ti$) or a small  maximum number of fixed-point iterations ($\maxit$) is chosen in the Algorithm \ref{decoupled_algo}.

\begin{algorithm}[htb!]
\DontPrintSemicolon
\caption{Decoupling algorithm: Inexact Newton scheme}
\label{decoupled_algo}
\SetAlgoNoLine
$n=0$\;
\Repeat {$\rn < \tn$} {
{$i=0$}\;
\Repeat{$\ri<\ti$ or $i=\maxit$}
{
{compute Newton updates $(\delta \du^{i+1},\delta \dv^{i+1})$ by solving \eqref{linearized decoupled problem_stat}}  \;
{evaluate the residual $\ri$ of the linear system \eqref{linearized problem_stat} }\;
{$i=i+1$}\;
}
{update the iterate $\du^{n}$ and $\dv^{n}$ by $\delta u^{n+1}$ and $\delta v^{n+1}$}\;
{ evaluate the  residual $\rn$ of the nonlinear system  \refAB} \;
$n=n+1$\; 
}

\end{algorithm}

This decoupling method is compared for different parameters in numerical tests of Section \ref{section:stat_results} to the fully coupling Newton method.
\label{chapter:multigrid}

\subsection{Multigrid scheme for the coupling approach}
\label{mg for coupled}
In this section we introduce a multilevel preconditioner which can cope with the strong coupling between PDE and ODEs. Such a coupling arises in the solver of the linear subsystems if the fully coupling Newton method is used instead of a decoupling scheme.  Coupled problems are commonly preconditioned by block preconditioning approaches, e.g.\  by simple block diagonal methods or a preconditioning of the Schur complement \cite{bockmandel}.  We will not use a block preconditioning approach because of the small dimension of the ODE part but instead we set up a coupling preconditioner based on the linear multigrid method.

{\color{cor}
  In fact, it is well known that the most efficient preconditioner for the PDE block is a multilevel preconditioner, because the number of linear iterations becomes independent of the mesh refinement \cite{Bramble:1993,Hackbusch:1985}. 
}
We consider a hierarchy of meshes \grd, where the index $0$ denotes the root mesh, i.e.\ the coarsest mesh from which all other meshes are derived by refinement. 
In this section we use the following notation for the system matrix of \eqref{linearized problem_stat}:
\begin{align}
\label{matrix K}
{\cal K}_l :=\begin{pmatrix} 
A_{h,u}^{\prime,l} & A_{h,v}^{\prime,l}\\
B_{h,u}^{\prime} & B_{h,v}^{\prime}
\end{pmatrix},\end{align}
where the index $l$ indicates the grid refinement level. 
The diagonal block $B_{h,v}^{\prime}$ does not depend on the mesh level, whereas the block $B_{h,u}^{\prime}$ does depend on the mesh level through the coupling term $\ut_h$ on the cell boundary.
Nevertheless, we do not use the notation with superscript $l$ in the blocks of the ODE part. 
In fact, we use an approximation in the coupling. 
To reduce the computational costs and simplify the implementation, the coupling ODE/PDE block is calculated at each level with the term $\ut_h$ computed at the finest level. In this way the whole ODE part does not depend on the refinement level $l$.
Our numerical results have indicated that this modification does not influence the performance of the multilevel algorithm.

The multigrid scheme is used as a preconditioner for a Krylov method applied to the system matrix ${\cal K}_l$. We use a  generalized minimal residual (GMRES)  method  because of the asymmetry of the system matrix, but a different Krylov method as, e.g.,\ the BiCG or BiCGStab would  also be appropriate for our purpose. 
In Section \ref{section:stat_results} we show numerically that the efficiency of the preconditioner is independent of the mesh size.
The work presented here is based for the PDE part on previous work by Janssen and Kanschat \cite{JanssenKanschat:2011}.
{\color{cor}
For the coupling PDE/ODE a new implementation was necessary. The difference to a standard multigrid implementation is that the transfer operators have to take into account the structure of the problem and treat the PDE and ODE parts in different way as explained below.
}
\subsubsection{Transfer operators}
For the transfer operators we use the following notation
\begin{align}
\label{transfer PDE}
R_l^{l-1} : V_l \rightarrow V_{l-1} ~({\rm restriction}), \quad P_l^{l-1} : V_{l-1} \rightarrow V_l ~({\rm prolongation}).
\end{align}
The restriction and prolongation operators act only on the PDE part, i.e.\ the finite element discretization, while the ODE part is transferred by the identity in both directions. 
The restriction and prolongation for the PDE part are implemented as intergrid transfers induced by the natural embedding of hierarchical meshes \cite{JanssenKanschat:2011}.
In matrix notation the restriction of the whole residual is given by application of the operator
\begin{align}
\label{transfer}
\begin{pmatrix}
R_l^{l-1} & 0\\
0 & I
\end{pmatrix},
\end{align}
where $I \in \mathbb R^{\NO\times \NO}$ denotes the identity matrix of the ODE part. The prolongation of the whole residual is defined analogously.
\subsubsection{Smoother}
\label{smoothers}
  In case of strong coupled problems, for example in fluid structure interaction, a common strategy for the smoothing process is to consider the full coupling only at the coarsest level and to smooth the two parts separately (decoupled) on the finer levels \cite{vanBrummelen:2008,Gee:2011,Langer:2014,Richter:2015}.
Since in our case the ODE part is small in comparison with the PDE one, we expect that the marginally more expensive smoothing of the whole coupled system at all levels would be efficient in the case of strong coupling. Therefore, we have compared the two strategies of (1) smoothing the whole system or (2) smoothing only the PDE part. For this comparison every efficient smoother would be appropriate, we have chosen the incomplete LU factorization (ILU). The two smoothers are denoted $S_1$ and $S_2$:
\begin{description}
\item $S_1$: Incomplete factorization (ILU) of the whole matrix ${\cal K}_l$ \eqref{matrix K};
\item $S_2$: incomplete factorization (ILU) of the PDE block as part of a Block Gauss Seidel scheme.
\end{description}
\subsection{Multigrid scheme for decoupling approach}
\label{mg for decoupled}
In the decoupled formulation as depicted in Figure \ref{fig:coupling} we use the multigrid only for the PDE part. Therefore, we use the standard multigrid scheme from \cite{JanssenKanschat:2011} with ILU as smoother.

\section{Numerical Results}
\label{sec:results}
In this section we make a comparison between the different numerical schemes presented in the previous section.  The following computations were performed using the C++ library deal.II, see Bangerth et al. \cite{dealII82}, with the UMFPACK library applied as direct solver on the coarse grid level \cite{davis2004algorithm}. Further implementation details can be found in \cite{gerecht}.

\paragraph{Intercellular signalling model}
\label{subsection:model}

Exemplarily we focus on a model for signaling of Interleukin-2 (IL-2) between T cells in the lymph node first presented by Busse et al. \cite{Busse:2010}. 
In the lymph node there are different types of cells. We consider in our test model IL-2 secreting cells and responder cells, i.e.\ cells that secrete IL-2 and cells that receive the IL-2 which diffuses in the extracellular region. The cells are fixed in space and on their surface there are receptors responsible for recognizing antigens.
A description of the modeled biological processes is given in Section \ref{application}. 
Here, we give a mathematical description of the model, which consists of a reaction-diffusion equation describing the distribution of IL-2 between the T cells in the intercellular area $\Omega$ coupled with ODEs for the intracellular processes involving the receptor dynamics. The coupling is described by a Robin boundary condition for each T cell indexed by $i=1,\dots,N_c$.

Let us consider the following notation:
\begin{itemize}
\item $u(t,x):[0,T]\times \Omega \to \R$ describes the concentration of IL-2 in the intercellular area.
\item $R_i(t),\BB_i(t)$ and $E_i(t) : [0,T]\to\R$ describe the number of IL-2 receptors, built receptor-complexes and internalized complexes for each of the simulated T cells. 
The receptors are distributed homogeneously on the cell surfaces.
\end{itemize}

\begin{subequations}
The mathematical model consists of one PDE
\label{problem}
  {\color{cor}
\begin{align} 
\begin{split}\label{pde_problem}
\partial_t u(t,x)&= \dd \Delta u(t,x)-\kd u(t,x) \quad \text{for all } (t,x) \text{ in } (0,T]\times\Omega,\\
D  \partial_n u(t,s)&=q_i -\kon R_i(t) u(t,s) +\koff \BB_i(t)\quad\text{for all $(t,s)$ in $(0,T]\times\Gamma_i$},\\
\partial_n u(t,s)&=0 \quad \text{for all $(t,s)$ in $(0,T]\times\Gamma_{out}$},
\end{split}
\end{align}
}
coupled with three ODEs for each T cell
\begin{align} 
\begin{split}
 \partial_t R_i(t)&=\vnoti+\vonei\frac{\BB_i(t)^3}{\kk^3+\BB_i(t)^3} 
 -\kono R_i(t) \widetilde{u}_i(t)\newln
-\kir R_i(t)+\koffo \BB_i(t)+\krec E_i(t)\quad\text{for all cells $i=1,\dots ,N_c,$}\\
 \partial_t \BB_i(t)&= \kono R_i(t) \widetilde{u}_i(t) - (\koffo+\kib)\BB_i(t),\\
 \partial_t E_i(t)&=\kib \BB_i(t)-(\krec+\kdeg)E_i(t)\label{bio_equation},\\
\ut_i(t)&=\frac{\int_{\Gamma_i}u(t,s)\text{ ds}}{\abs{\Gamma_i}}, 
\end{split}
\end{align}
\end{subequations}
with given initial conditions for $u(0),R_i(0),\BB_i(0)$ and $ E_i(0)$ for all cells  $i=1,\dots,N_c$.  The used parameters and their values are described in Table \ref{table_parameters}. 
The parameters and the chosen initial values have been taken from \cite{Busse:2010}.
\begin{table}[htb!]
   \centering
 \caption{Model parameters} \label{table_parameters}
    \begin{tabular}{@{}lll@{}}
       \toprule  
Symbol&Value&Parameter\\
\midrule
$q_i$&$0-22000$ mol./cell/$h$&IL-2 secretion rate\\
$\dd$&36000 $\mu m^2/h$&{Diffusion coefficient of IL-2}\\
$\kd$&0,1/$h$&{Extracellular IL-2 degradation}\\
$\vnoti$&150 mol./cell/$h$&{Antigen stimulated IL-2 receptor expression rate}\\
$\vonei$&3000 mol./cell/$h$&{Feedback induced IL-2 receptor expression rate}\\
$\kk$&1000 mol./cell&{Half-saturation constant of feedback expression }\\
$\kon$&111,6 /nM/$h$&{IL-2 association rate constant to IL-2 receptors}\\
$\koff$&0,83/$h$&{IL-2 dissociation rate constant from IL-2 receptors}\\
$\kir$&0,64/$h$&{Internalization rate constant of IL-2 receptors}\\
$\kic$&1,7/$h$&{Internalization rate constant of receptor complexes}\\
$\krec$&9/$h$&{Recycling rate constant of IL-2 receptors}\\
$\kdeg$&5/$h$&{Endosomal degradation constant IL-2 receptors}\\
$r$&5$\mu m$&{Cell radius}\\
$d$&5$\mu m$&{Cell to cell distance }\\
\bottomrule
 \end{tabular}
\end{table}
The two types of cell that we consider share the same receptor dynamics but differ in the IL-2 secretion rate:
\begin{itemize}
  \item Secreting T cells, which {\color{cor}secrete} IL-2 with the secretion rate $q_i=2500 \text{ mol}/h$,
\item Responding T cells with $q_i=0$.
\end{itemize}
\begin{table}[ht]
    \centering
    \caption{Value of the IL-2 degradation parameter and corresponding maximum eigenvalue of $J$, see \eqref{Jacobian}}
    \label{tb:parameters}
    \begin{tabular}{ccc}
      \toprule  
      &  $k_d$ & $\lambda_{max}$\\
    \hline
   Biologically relevant test model & 0.1 & 8.8\\
    \hline
   Artificial test model& 1000 & 0.01\\
    \bottomrule
    \end{tabular}
\end{table}
In Section \ref{application} we show a simulation for a microenvironment configuration with 216 T cells, among them 54 randomly chosen secreting T cells. That simulation is used to give a biological interpretation to the results obtained with the presented model.
The numerical tests shown here are performed using a subproblem. We consider only eight cells among them only one is secreting IL-2. The reduction in size allows for fast tests and does not affect the study of the strength of the coupling. In fact, we have observed by calculations using many different configurations that the strength of the coupling depends on the model parameters and not on the number of cells.
The configuration is displayed in Figure \ref{3d_pic2} where the responding T cells are in grey and the secreting T cell is highlighted.
  For this test problem we found a stationary state numerically. 

  \medskip
{\color{cor}
The goal of the next subsections is to show the numerical comparison between the two schemes depicted in Figure \ref{fig:coupling}.
We have divided the results in three subsections:
\begin{itemize}
  \item Subsection \ref{preconditioner_results}: at first we compare the two smoothers $S_1$ and $S_2$ that we have chosen for the coupling solver, see again Figure \ref{fig:coupling}. 
\end{itemize}
Once we have decided which is the best smoother for the coupling solver, we compare this solver with the decoupling solver in two cases: 
\begin{itemize}
  \item Subsection \ref{section:stat_results}: the stationary case;
  \item Subsection \ref{section:instat}: the non-stationary case.
\end{itemize}
Since the results depend on the strength of the coupling, we want to test our solvers for \emph{strong} and \emph{weak} coupling conditions for both the stationary and the non-stationary cases. In the stationary case, we consider two configurations changing the parameter $k_d$. The two values of $k_d$ and the respective maximal eigenvalues of the Jacobian are shown in Table \ref{tb:parameters}. 
We denote the case with low value of $k_d$ as \emph{biologically relevant} since it is in the range of reference values given in Table \ref{table_parameters}. 
According to the criterion \eqref{strength criterion strong} it is a strongly coupled system.

We refer to the case with a high value of $k_d$ as an \emph{artificial} setting since this value is far larger than the biological assumptions taken in other publications with this model \cite{Busse:2010,thurley14}.
According to the criterion \eqref{strength criterion weak} it is a weakly coupled system.
High values of $k_d$ correspond to large degradation of $u$.
In this way the influence of $v$ on $u$ is diminished.
Thus the PDE part is artificially decoupled from the ODE part.
We use the two proposed settings to compare the solvers.
In the next paragraphs we give a biological interpretation of the {\color{cor1}two settings}.

\paragraph{{\color{cor1}
Biological scenario for the strongly coupled (biological relevant) test model}}
{\color{cor1}
  This model consists of 8 T cells among which one is a secreting cell. All parameters are like in the real application, in the lymph node (Table \ref{table_parameters}). After the definition in \eqref{strength criterion strong} we have a strongly coupling between the equations. In Figure \ref{x} (a) we see the time course of the IL-2R receptors (the sum of R and C receptors) on the surface of each cell. The receptors need some minutes to get the IL-2 for a slight increase and reach their steady state quite early after 9 hours. In this configuration one secreting cell is not enough to activate the others.  In Figure \ref{x} (b) we see the time course of the averaged IL-2. Again we have a slight increase in the first 3 hours and the steady state regulates quite early after 9 hours. The IL-2 concentration is low, around 0.045 pM.
}
\begin{figure}[H]
  \caption{Dynamic behavior for the biological relevant test model which is strongly coupled (8 T cells among them one secreting)}
\centering
\subfigure[Time course of the amount of IL-2R receptors (R+C) on the surface of each cell. No cell is activated in this configuration]{
  \includegraphics[width=.7\textwidth]{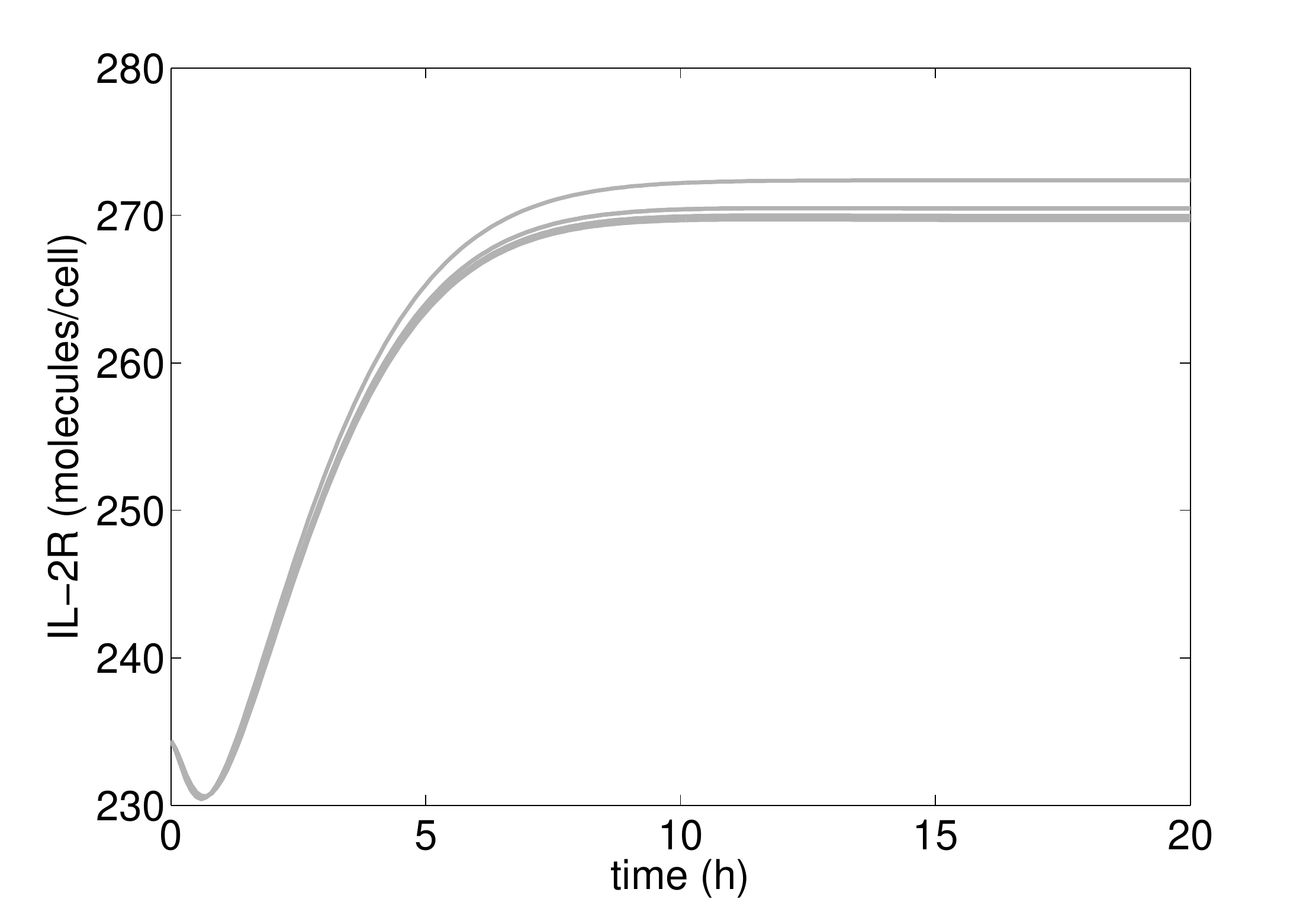}}
  \subfigure[Time course of the averaged IL-2 concentration on the surface of all cells]{
  \includegraphics[width=.7\textwidth]{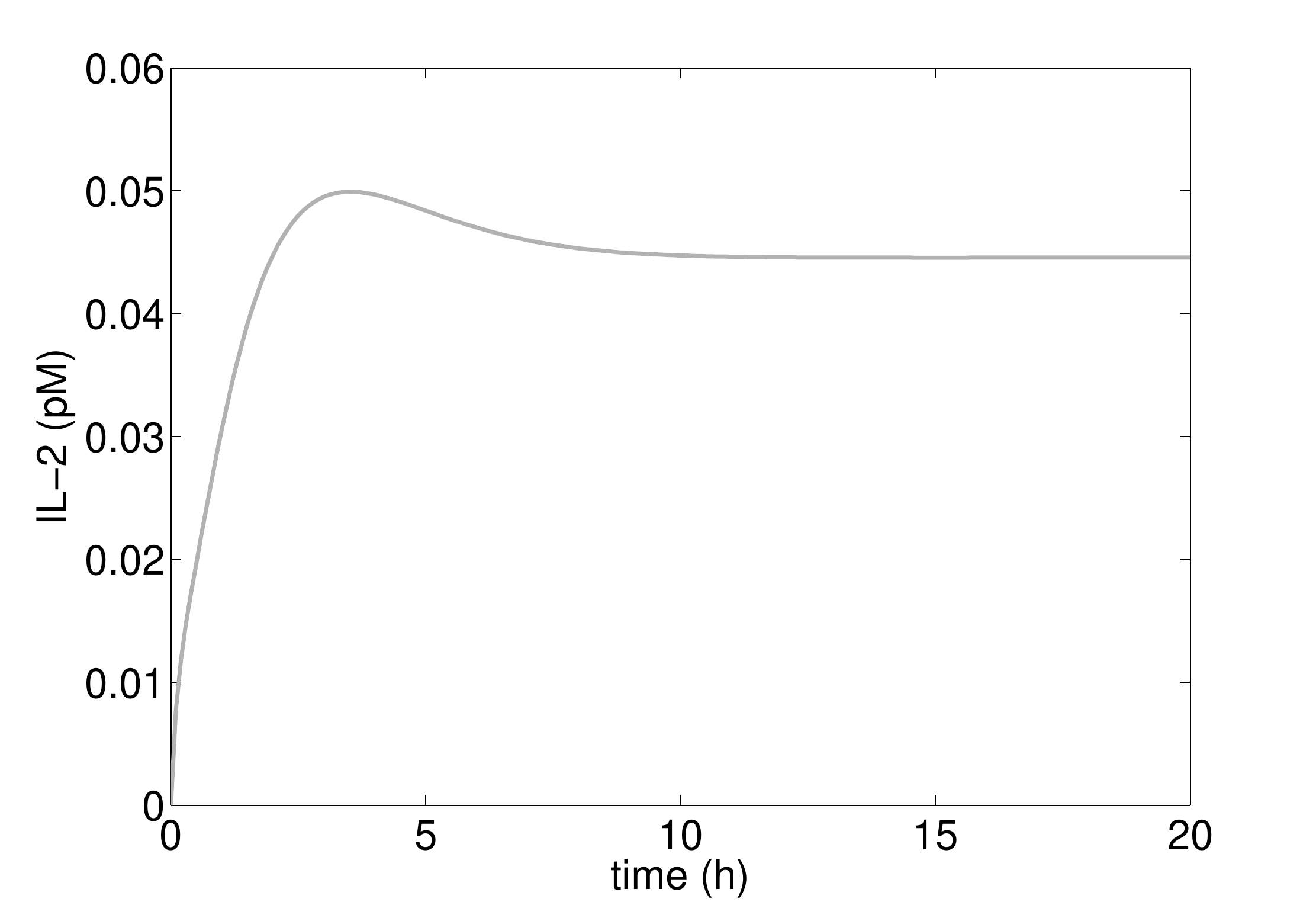}}
  \label{x}
\end{figure}

\paragraph{{\color{cor1}Biological scenario for the weakly coupled (artificial) test model}}

This model differs from the biological relevant test model only in one parameter, $k_d$ , the degradation rate.
The parameter is chosen artificially to create a weekly coupling between the equations.
It is chosen to be greater by four orders of magnitude with respect to the biologically relevant value (Table \ref{tb:parameters}). Therefore, IL-2 degrades $10^4$ faster than in the strongly coupled scenario. 
It follows that there is much less IL-2 available for the responder cells than in the strongly coupled model. Numerical simulations ({\color{cor1} Figure \ref{x+1}}) show that the amount of IL-2 is reduced by three orders of magnitude. 
{\color{cor1} 
  Thus, there will be less IL-2 for the responder cells and also no chance for activation. In the steady state which occurs immediately, whereas in the strongly coupled model it occurs after 9 hours, there is only $2.2\cdot 10^{-4} \text{~pM}$ IL-2 which is three orders of magnitude less than in the biological relevant scenario.
The total amount of IL-2R receptors stays almost constant for all times.
}
With this concentration changes in the ODEs the secreting cell also produces less IL-2.

In the non-stationary case, we do not have to change $k_d$ to compare two configurations. Instead we use two different time discretizations with steps $\Delta t = 0.01 \rm h$ and $\Delta t = 0.1 \rm h$. By reducing the time step in this way we also reduce the strength of the coupling as shown in the results.

\begin{figure}[H]
  \caption{Dynamic behavior for the artificial test model which is weekly coupled (8 T cells among them one secreting)}
\centering
\subfigure[Time course of the amount of IL-2R receptors (R+C) on the surface of each cell. No cell is activated in this configuration]{
  \includegraphics[width=.7\textwidth]{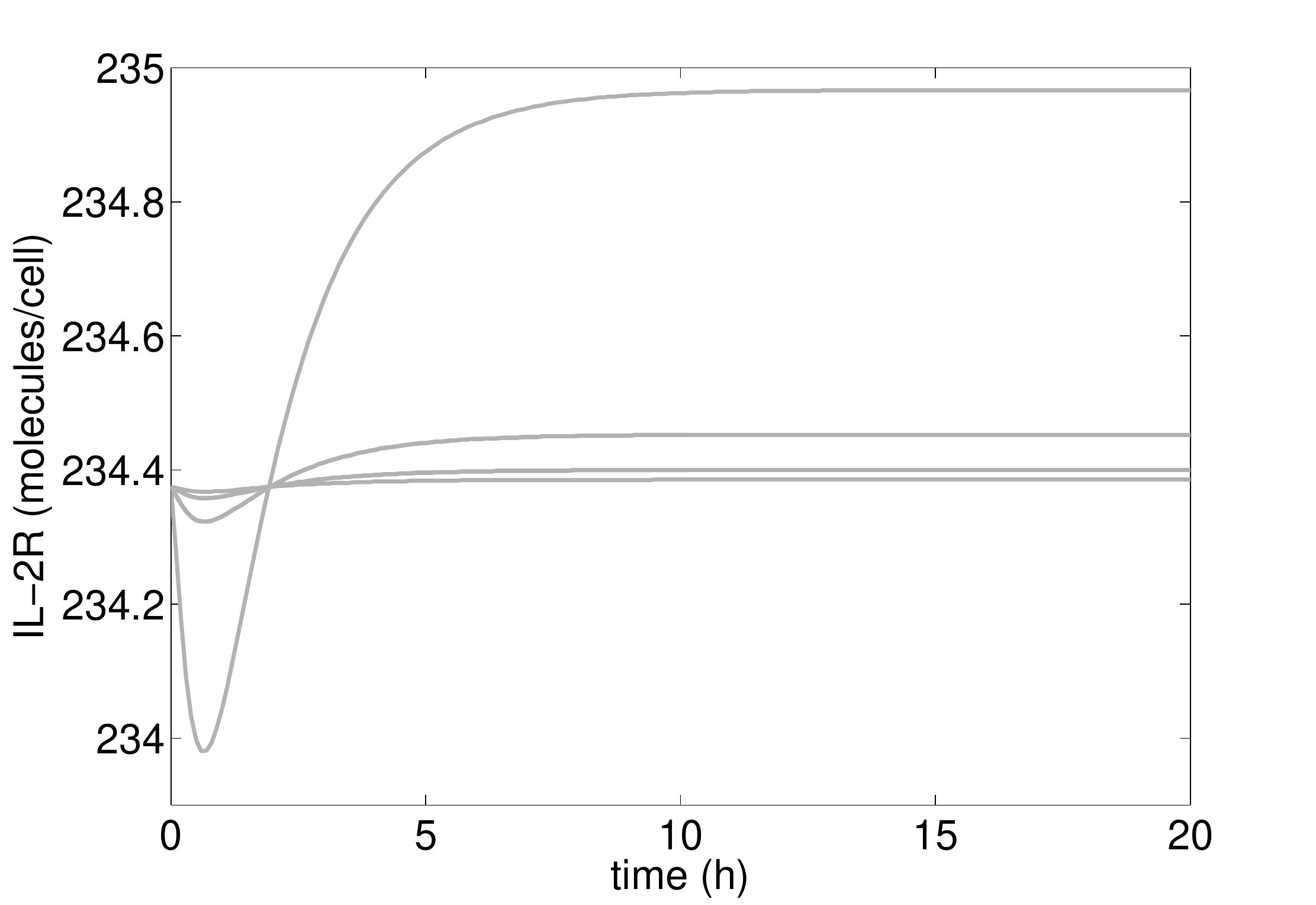}}
  \subfigure[Time course of the averaged IL-2 concentration on the surface of all cells]{
  \includegraphics[width=.7\textwidth]{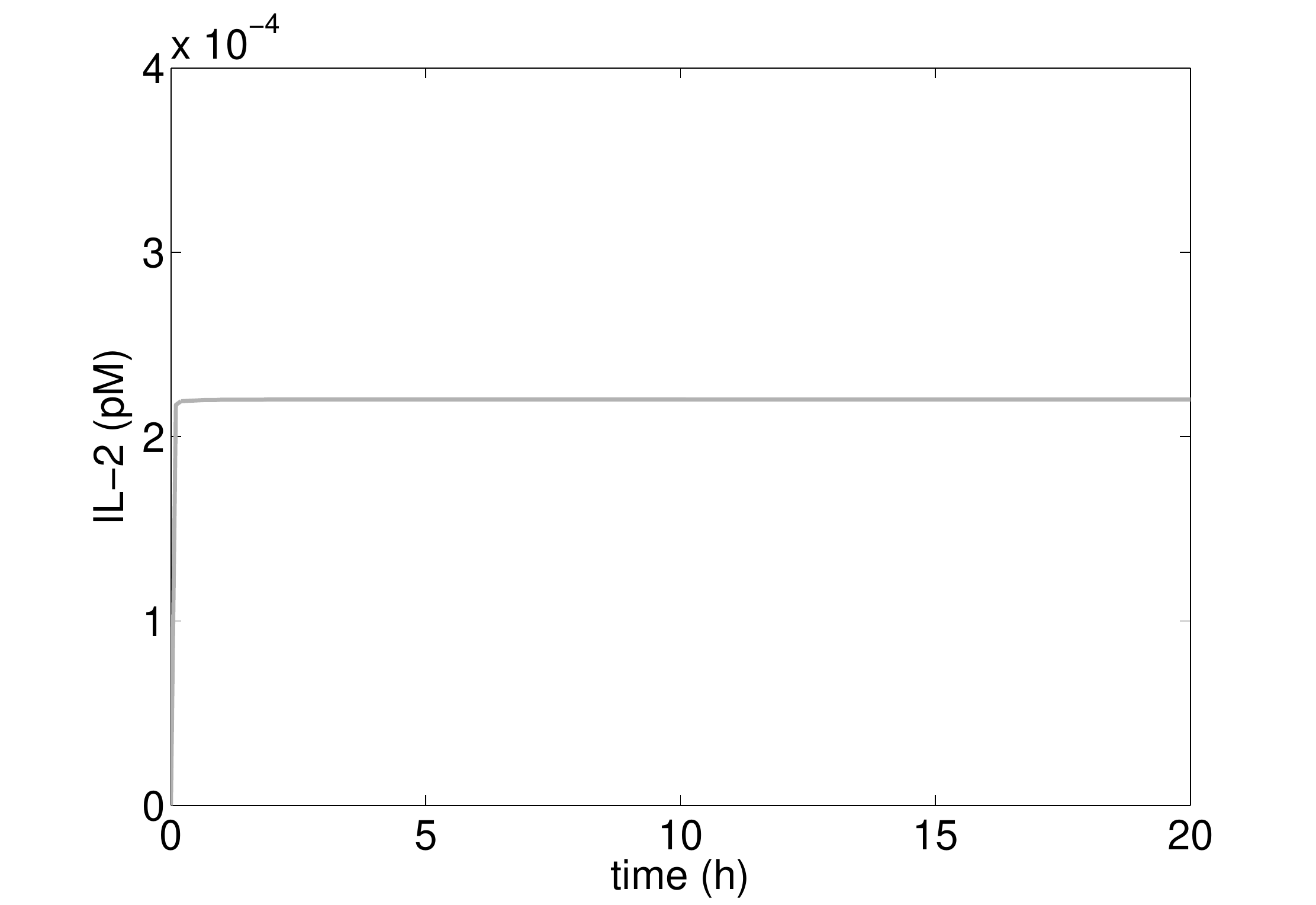}}
  \label{x+1}
\end{figure}
}

\subsection{Choice of multigrid preconditioner for the coupling solver}
\label{preconditioner_results}
The two smoothers have been described in Subsection \ref{smoothers}. We compare them using the configuration with a strong coupling strength, i.e.\ the biologically relevant case. Recall that $S1$ uses a factorization of the whole matrix, while $S2$ uses a factorization of only the PDE part.

We compute the number of GMRES steps over all {\color{cor}Newton} steps ($\Sigma{n}$)  and the average reduction rate ($r$) of the residual in each GMRES step.
We have set the Newton accuracy to $10^{-6}$.
The accuracy of the linear solver {\color{cor1} per Newton step} is set to $10^{-11}$.
We have observed that the number of Newton steps only depends on the coupling, the nonlinearity of the equation and the accuracy of the solver. In particular it is not dependent on the grid refinement. The latter effect is known as ``asymptotic mesh independence'' and a proof of it can be found in \cite{Weiser:2005,Allgower:1987} for Newton-type methods applied to nonlinear operator equations and, e.g., in \cite{Heinkenschloss:1993} for a Gauss-Newton method for nonlinear Least Squares problems.
The number of smoothing cycles are set in both smoothers to three. We have observed that a larger number of smoothing steps is unnecessary to improve the smoother performance.
To test the mesh independence of the multigrid preconditioner we compute several refinement levels up to five refinements as indicated in the first column of Table \ref{tb:precon}.
We use a global refinement of the grid. The finest grid has $885673$ degrees of freedom, the coarsest 367 degrees of freedom for the PDE part. Additionally 24 degrees of freedom for the ODE part are coupled to the PDE part.
%
Table \ref{tb:precon} shows that the number of GMRES steps preconditioned by the two smoothers is almost constant for $S_1$ (from 46 at $L=2$ to 54 at $L=5$) and increases moderately for $S_2$ (from 69 at $L=2$ to 84 at $L=5$). Therefore we have chosen $S_1$ because it is mesh independent and more effective. In fact, the computational costs per smoothing iteration are almost the same for the two smoothers because the ODE block is not large. Since the costs are comparable, we see that the smoother $S_1$ costs from 33\% (at $L=2$) to 35\% (at $L=5$) less than the smoother $S_2$.
\begin{table}[ht]
    \centering
 \caption{Reduction rates of different preconditioners }
    \label{tb:precon}
    \begin{tabular}{@{}cccccc@{}}
      \toprule  
&\multicolumn{2}{c}{MG-$S_1$ (ILU)}&\multicolumn{2}{c}{MG-$S_2$ (Bl.-ILU)}
\\
\cmidrule(lr) {2-3}\cmidrule(lr){4-5}

$L$& $\log_{10}r$  
&$\sums$&$log_{10}r$&  $\sums$ \\
\hline   
2
  &2.00 &46& 1.41&69\\
3
  &1.92 &51 & 1.27&77\\
4
 &1.85 &54& 1.21&81 \\
5
 &1.81 &54&  1.15&84  \\
\bottomrule
\multicolumn{6}{l}{{\it Notation:}  $\sums$ GMRES iterations}\\
\multicolumn{6}{l}{{\bf $\quad$$\quad$ $\quad$$\;\;$$\quad$}  in all Newton steps}\\
\multicolumn{6}{l}{{\bf $\quad$$\quad$ $\quad$$\;\;$}  $r$ average reduction rate}\\
\multicolumn{6}{l}{{\bf $\quad$$\quad$ $\quad$$\;\;$}  $L$ refinement level }
 \end{tabular}    
\end{table}
\subsection{Stationary case}\label{section:stat_results}
In this section we compare the two approaches depicted in Figure \ref{fig:coupling} where for the coupling method (a) we use the smoother $S_1$ as explained above.
For the coupling approach we use an iterative method preconditioned by the multigrid method described in \ref{mg for coupled}, while for the decoupling approach we use the method described in \ref{mg for decoupled}.
In the coupled case, since the matrix is asymmetric we use a GMRES method.
In the decoupling approach the system matrix is symmetric, therefore the solver of choice would be, e.g.\ the conjugate gradient method (CG). Nevertheless, for a direct comparison we use instead the GMRES method also for this approach. In fact, we have observed that, in combination with the preconditioner, both solvers have similar performance.
Furthermore, to make the schemes comparable we use the same accuracy of the GMRES solver set to $\ti$ for both schemes. In this way the number of Newton steps to solve the nonlinear problem is independent of the approach and we can compare the total number of GMRES steps to solve for a Newton accuracy of $\tn$.

In Table \ref{tb:coupling} we compare the number of  Newton steps ($n$) and the number of total \mbox{GMRES} iterations ($\sums$) needed to obtain a solution of accuracy ($\tn=10^{-6}$). In each Newton step the decoupling scheme described in Algorithm \ref{decoupled_algo} is iterated until a residual $\ri < \ti=10^{-11}$ is reached. The parameter $\maxit$ is not used in this test, we will use it later in the non-stationary case.
The average \mbox{GMRES} iterations per Newton step is denoted by $\bar{s}$ and the sum over all \mbox{GMRES} steps by $\sums$. We globally refine the coarse grid three times up to a number of $114929$ degrees of freedom. 
\begin{table}[ht]
    \centering
  \caption{Coupling vs decoupling solver}

    \label{tb:coupling}
    \begin{tabular}{@{}rcccccc|cccccc@{}}
      \toprule  
&\multicolumn{6}{c|}{``biologically relevant'' problem $\lambda=8.8$}&\multicolumn{6}{c}{``artificial'' problem  $\lambda=0.01$}\\\midrule
&\multicolumn{3}{c}{decoupling}&\multicolumn{3}{c|}{coupling}&\multicolumn{3}{c}{decoupling}&\multicolumn{3}{c}{coupling}\\
$L$&  $n$&$\bar{s}$& $\sums$    &$n$&$\bar{s}$& $\sums$  &$n$&$\bar{s}$& $\sums$    &$n$&$\bar{s}$& $\sums$  \\\midrule   
 2  &7 &748 &5236 &7& 6.6 & 46& 3&11 & 33&3 &7 &21\\
3  &7 &921 &6444 &7& 7.3 & 51& 3&11 & 33&3 &7 &21\\
 4&7 &957    &6699     &7& 7.7 & 54& 3&11.7 & 35&3 &7 &21\\
\bottomrule
 \multicolumn{13}{l}{{\it Notation:}  $\sums$ GMRES iterations during all Newton steps}\\
\multicolumn{13}{l}{{\bf $\quad$$\quad$ $\quad$$\;\;$}  $n$ Newton steps}\\
\multicolumn{13}{l}{{\bf $\quad$$\quad$ $\quad$$\;\;$}  $\bar{s}$ average GMRES iterations per Newton step}\\
\multicolumn{13}{l}{{\bf $\quad$$\quad$ $\quad$$\;\;$}  $L$ refinement level }\\
\multicolumn{13}{l}{{\bf $\quad$$\quad$ $\quad$$\;\;$}  $\lambda$ largest eigenvalue of the sensitivity matrices }
\end{tabular}   \end{table}
It can be observed that the coupling solver outperforms the decoupling one in both cases, the ``biologically relevant'' and the ``artificial'' configuration. While in the \emph{strongly} coupled case the coupling solver is much better than the decoupling approach (compare the $4^{\rm th}$ column with the $7^{\rm th}$ column), in the \emph{weakly} coupled case it leads to 30--40\% less GMRES iterations (compare the $10^{\rm th}$ column with the $13^{\rm th}$ column).

\subsection{Non-stationary case}\label{section:instat}
In this section we solve the time-dependent system \eqref{eq:instationary_problem} until a stationary solution is reached. In each time step a coupled non-linear system has to be solved. In contrast to the stationary case, the strength of the coupling  is reduced when using small time steps. 
For the decoupling approach it is less effective (in the non-stationary case) to solve the linear problems up to a very small tolerance ($\ti$) as can be seen from the numerical results. For this reason we have introduced the parameter $\maxit$ in Algorithm \ref{decoupled_algo}, and for comparison purposes we have made numerical tests with the values $\maxit =1,\dots,4$, see first column in Table \ref{tb:instat_coupling}.
\begin{table}[ht]
    \centering
  \caption{Decoupling and coupling solving of the non-stationary problem }
    \label{tb:instat_coupling}
    \begin{tabular}{@{}lccc|clc@{}}
      \toprule  
      &&\multicolumn{2}{c}{$\Delta t=0.1 \rm h$}&\multicolumn{2}{c}{$\Delta t=0.01 \rm h$}&$$\\
&$\maxit$&  $\sumn$& $\sums$ &  $\sumn$& $\sums$&$$\\
\midrule
decoupling& 1&1393 &3375& 5566&14016 \\
           &2&754 &3792&3395&15217 \\
           &3&547 &4363&2880 &21153\\
           &4&448 &4775&2871 &23658\\
\hline   
coupling&&356&1799&2868&11299\\
\bottomrule
\multicolumn{7}{l}{{\it Notation:}  $\sumn$ Newton steps in all time steps}\\
 \multicolumn{7}{l}{{\bf $\quad$$\quad$ $\quad$$\;$}  $\sums$ Krylow iterations in all Newton steps}\\
 \multicolumn{7}{l}{{\bf $\quad$$\quad$ $\quad$$\;$}  $\maxit$ maximum of iterations}\\
 \multicolumn{7}{l}{{\bf $\quad\;\;$ $\quad\quad\;\;\quad\quad\quad$$\;$}  per Newton step}
 \end{tabular}   
\end{table}
This table reports as well the sum of computed newton steps ($\sumn$) and the sum of computed GMRES steps ($\sums$) over all time steps. 
The results are listed for computations on a spatial grid level $L=2$ (2189 degrees of freedom) with 200 ($\Delta t=0.1 \rm h$) or 2000 ($\Delta t=0.01 \rm h$) time steps, with final time 20 hours.

A higher maximal number of fixed-point iterations per Newton step (i.e. a higher $\maxit$) increases the accuracy of the linear solver and thus reduces the  number of Newton steps. The decoupling solving scheme with a $\maxit=4$ results in a number of Newton steps not much larger than in coupling solving scheme (448 vs 356) but with more than twice the number of computed GMRES steps (4775 vs 1799).
Furthermore, it can be observed that the number of total GMRES steps decreases reducing $\maxit$, see the $3^{\rm rd}$ column.
A comparison of the computational time should consider also the time per Newton step. In fact, we observe that the computing time for solving with the decoupling approach is much larger than the double of the time needed by the coupling approach. In fact, associated to each Newton iteration there are additional computational costs, e.g.\ building or updating the Jacobian and the residual. Since the time per Newton step depends on the specific implementation (e.g.\ how often the Jacobian is updated), we restrict the comparison to the total number of linear solving. For this reason we consider a better choice the value of $\maxit=1$.

As already remarked, the effectiveness of the decoupling solver depends on the strength of the coupling and thus on the size of the time step.
In fact the coupling solver needs for time steps $\Delta t=0.1 \rm h$ around half of the iterations of the decoupling solver, while for smaller time steps ($\Delta t=0.01 \rm h$) the iterations of the coupling solver are reduced by $20 \%$  compared to the decoupling solver with $\maxit=1$. 
\begin{remark}\label{remark:CN}
By using the implicit Euler scheme we have shown that the coupling solver is more effective.
The use of a higher order time scheme, e.g.\ the Crank-Nicolson scheme, allows for larger time steps to produce the same accuracy. As shown above, larger time steps lead to a stronger coupling during the time integration. We expect therefore that the coupling solver is even more effective using a higher order time scheme.
\end{remark}
\begin{remark}
For the test that we have done in this work we did not use a globalization method for the Newton convergence. Our experience is that for the dynamic problem we do not need globalization, while for the stationary problem some critical configurations could be solved only using pseudo time steps approach to get a good start solution for Newton to converge.
\end{remark}

\paragraph{Solver performance with increasing number of T cells}
To test how the solver scales with the number of T cells, i.e.\ how the total number of Newton steps and the total number of GMRES steps increase with the increasing number of the simulated T cells, we have run simulations with 8, 27 and 64 T cells. The parameter set is the same as that for the previous simulations, in particular we set the final time to 20 hours and a time step $\Delta t = 0.1 \rm h$. To maintain the same ratio between secreting and responding T cells, around 1/8 of the total number of T cells is randomly chosen as secreting, i.e.\ 1, 3, and  8 for the simulations with 8, 27 and 64 T cells respectively. Since the computing time for the decoupling solver is much larger than that for the coupling solver, we have run the simulations at refinement level $L=1$.

In Table \ref{tb:instat scaling} we report the total number of Newton steps and the total number of GMRES steps for all time steps. For the coupled problem we use the notation $\sumn_C$ and $\sums_C$ and for the decoupling solver the notation $\sumn_D$ and $\sums_D$. The last column shows the reduction of total GMRES steps by using the coupling solver with respect to the decoupling one.
It can be observed that the gain of the coupling solver scales perfectly with the number of T cells, i.e.\ the reduction of total number of GMRES steps by using the coupling solver is around 50\% independently of the number of T cells.

\begin{table}[ht]
    \centering
    \caption{Scaling with the number of T cells}
    \label{tb:instat scaling}
    \begin{tabular}{cccccc}
      \toprule  
& \multicolumn{2}{c}{Coupling} & \multicolumn{2}{c}{Decoupling} & Reduction of $\sums$ \\
    Num. T cells    &     $\sumn_C$ & $\sums_C$ & $\sumn_D$ & $\sums_D$ & $\sums_C / \sums_D$\\
\midrule
8 & 337 & 1601 & 1250 & 2997 & 53\%\\
\hline   
27 & 334 & 2007 & 1214 & 3389 & 59\%\\
\hline   
64 & 344 & 2316 & 1343 & 4386 & 52\%\\
\hline   
{\color{cor1}125} & {\color{cor1}345} & {\color{cor1}2672} & {\color{cor1}1325} & {\color{cor1}4813} & {\color{cor1}55\%}\\
\bottomrule
\multicolumn{6}{l}{{\it Notation:}  $\sumn$ Newton steps in all time steps}\\
 \multicolumn{6}{l}{{\bf $\quad$$\quad$ $\quad$$\;$}  $\sums$ Krylow iterations in all Newton steps}\\
 \end{tabular}   
\end{table}
\section{Application to a cellular microenviroenment}
\label{application}
We present a numerical result applied to a cellular microenvironment consisting of 216 T cells. This in silico experiment can be considered as representative to observe the cell-to-cell interactions in the lymph node. In our configuration 54 T cells are randomly chosen as T helper cells and release IL-2 locally into the immunological synapse with a rate of $q_i = 3500 ~\text{mol/h}$. The {\color{cor1}remaining cells} are responding T cells with $q_i = 0$ which absorb the IL-2 from the environment. 
The amount of IL-2 which is absorbed by each cell is proportional to the number of the expressed receptors on their surface and this number depends again on the amount of absorbed IL-2. 
This means, the more IL-2 a cell can absorb the more receptors are expressed on its surface.
If a certain amount of receptors are expressed, the T cell is activated and ready for proliferation and differentiation determining the type of the immune response. The more IL-2 a cell absorbs the better is the chance for it to get activated. Thus, the result is a competitive behavior between the cells for IL-2. 

\begin{figure}[H]
  \caption{Dynamic behavior of 216 T cells among them 54 secreting T cells}
\centering
\subfigure[Time course of the amount of IL-2R receptors (R+C) on the surface of each of the activated (yellow curves) and non-activated (gray curves) cells]{
  \includegraphics[width=.7\textwidth]{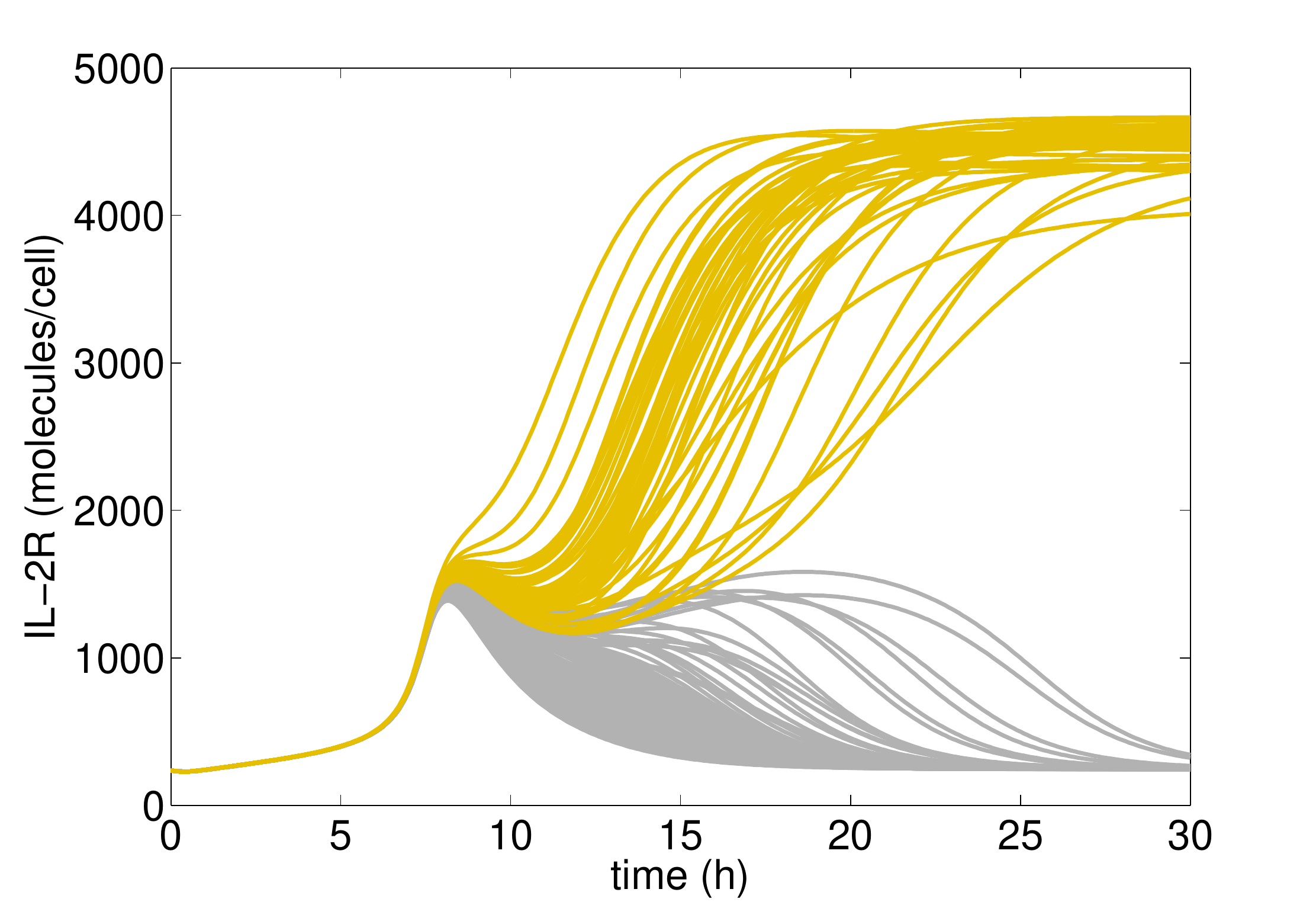}}
  \subfigure[Time course of the averaged IL-2 concentration on the surface of all activated (yellow curve; at t=9h is the curve above) and non-activated (gray curves) cell surfaces]{
  \includegraphics[width=.7\textwidth]{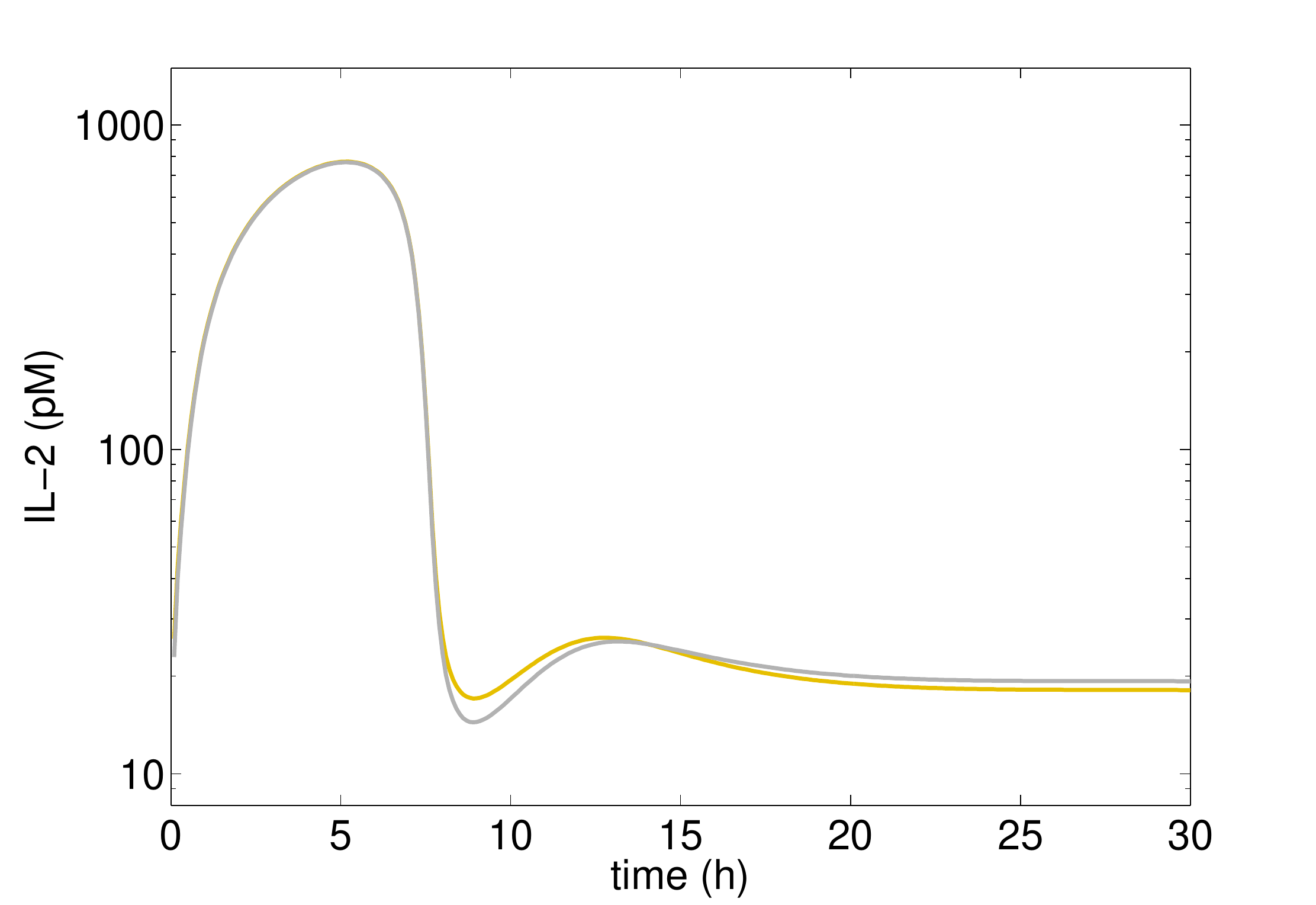}}
\label{IL-2 and IL-R plot}
\end{figure}

\begin{figure}[H]
  \caption{IL-2 concentration distribution in a cell microenvironment consisting of 216 T cells. This inhomogeneous IL-2 pattern develops in the steady state after 30 hours of activation of the IL-2 signaling pathway}
\centering
  \includegraphics[width=1.0\textwidth]{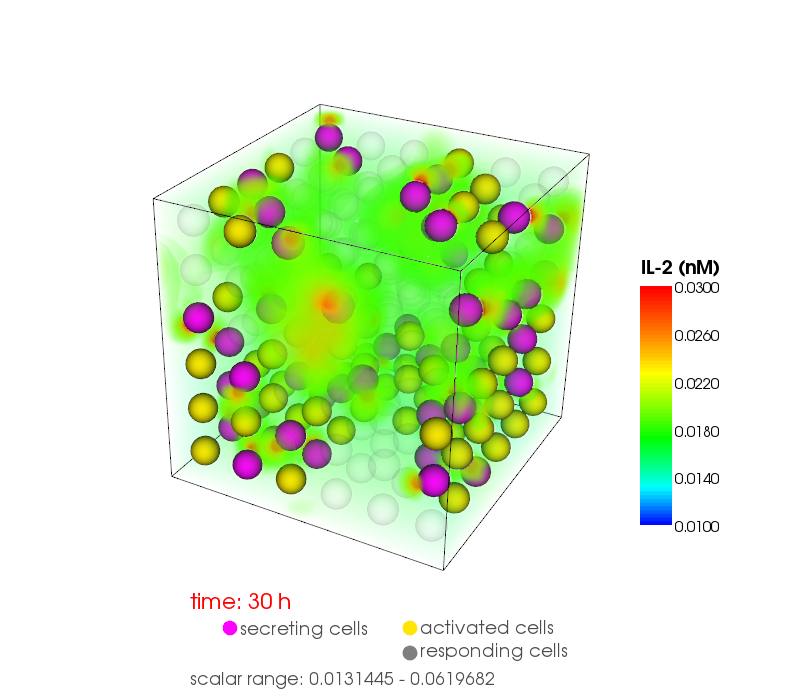}
\label{3D Visualization}
\end{figure}

For this {\color{cor1}application} we have used a final time of 30 hours, a time step $\Delta t = 0.1\rm h$ and a refinement level $L = 2$. {\color{cor1}Due to the value of $k_d = 0.1/h$ the equations are strongly coupled.} 
  In Figure \ref{IL-2 and IL-R plot} (a) we see the time development of the amount of receptors for each T cell. At  the beginning the amount of receptors grows {\color{cor1}strongly} in all cells with the same speed. 
  {\color{cor1}After 9 hours they have five times more receptors than in the biological relevant test model. Thereafter they split into two groups, the activated (yellow curves) and non-activated (grey curves) cells. A cell is called activated when it has more than 4000 IL-2R molecules.
  }
After 30 hours the T cells have either become activated by building new receptors or their number of receptors has decreased to the starting level. 
{\color{cor1} These biological processes take some time such that the steady state will be reached later after 30 hours.}
The Figure \ref{IL-2 and IL-R plot} (b) shows the averaged IL-2 concentration at the cell surface. Here, the same colors are used, yellow (see electronic version of the paper for the colors) for the activated and gray for the non-activated cells. 
Interestingly, there is almost no difference in the averaged amount of IL-2 between activated and non-activated cells except at the time point of decision, i.e. after 9 hours of activation. The cells which get slightly more IL-2 at this time will be activated. All others will down-regulate their receptors. 
{\color{cor1}In the steady state there is over 100 times more IL-2 than in the corresponding test model with one secreting cell and in the initial phase $10^4$ times more. One of the most important biological findings from these simulations was the heterogeneity of the IL-2 amount in time and space.}
Despite the rapid diffusion of IL-2, spatial inhomogeneities occur in the concentration distribution (Figure \ref{3D Visualization}) and large gradients develop over several orders of magnitude (Figure \ref{IL-2 and IL-R plot} (b)). In Figure \ref{3D Visualization} we present a volume visualization of our numerical results for the chosen configuration in the steady state performed {\color{cor1}with Covise \cite{Covise}}. 
Transparency of the whole data set is used to make the inner structure visible and a flexible mapping of the data on colors and opacity to visualize the different T cells for a realistic representation. The (randomly) chosen secreting cells are marked with magenta (see color version of the paper). By visualizing only a part of the IL-2 range we can see the secretion points on the cell surfaces (small red domains). There, we find the largest IL-2 concentration. The yellow cells are the activated cells and the transparent gray ones the responder cells which have down-regulated their receptors. Further visualization and numerical results are presented in \cite{thurley14}.

\section{Conclusions}
\label{conclusions}
In this paper, we considered a coupled nonlinear system consisting of a parabolic partial differential equation and many ordinary differential equations, which emerges e.g.\ in systems biology by modeling intercellular signaling pathways.
We presented numerical results for an application in immunology: the dynamics of cytokine (Interleukin-2) signaling between different types of T cells.
The presented methods can nevertheless be used for solving other signaling pathways or other applications modeled by systems of coupled PDE/ODE equations.

{\color{cor}Previous works dealing with coupled PDE/ODE models have presented mainly decoupling approaches (especially for biological applications) and have not shown a quantitative comparison with coupling approaches (see references in the introduction).
This paper showed in a systematic way a \emph{quantitative} comparison between coupling and decoupling approaches for this class of problems. 
}
Specifically, we used a sensitivity analysis and its numerical implementation to study the coupling/decoupling strategy. This approach, applied to the model for Interleukin-2 signaling, indicated that a coupling strategy performs better for our strongly coupled model.
We implemented a solution method based on a Newton-type solver with a multigrid preconditioner.
Depending on the time step length, the total number of linear solver iterations over all Newton steps can be reduced with this strategy up to 50\%.
The saving in computational time is typically much more than 50\% because also the number of Newton iterations is higher in the decoupling scheme in comparison with the coupling one, and associated to each Newton iteration there are additional computational costs. Nevertheless, the quantitative estimation of the computing time highly depends on how the Newton method is implemented. Therefore, in this work we have restricted the comparison to the total number of linear solving over all Newton steps.

{\color{cor1}As future work}, we indicate a possible strategy for an additional reduction of the computation time by using  local mesh refinement both in space and time. In particular, different time grids for the PDE and the ODE part allow, depending on the strength of the coupling, to decrease the number of time steps for the computationally expensive PDE part of the model. The essential question is how to choose the two time grids without decreasing the overall convergence rate of the method.
An a posteriori error estimator for the errors of the PDE and ODE discretization is necessary to reach a wanted accuracy efficiently by iterative adaptive refinement.
The use of a refinement strategy based on such an error estimator allows to control  the two time grids separately and obtain an optimal time discretization for both parts of the system.
The complex realization of such a method is subject of our current research.\\

\noindent \textbf{Acknowledgements}

The authors thank Prof.\ Rolf Rannacher for the numerous fruitful discussions on numerical aspects of this work. Furthermore, the authors acknowledge Prof.\ Thomas H\"ofer for his support and for provisioning this concrete application in immunology.
{\color{cor1}In addition, the authors gratefully acknowledged the work for the visualization in 3D done by Marcus Schaber (Visualization and Numerical Geometry Group at the Interdisciplinary Center of Scientific Computing (IWR), Heidelberg.)}\\ 
T.C. was supported by Deutsche Forschungsgemeinschaft (DFG) through the project CA 633/2-1.\\
E.F. and D.G. were supported by the Helmholtz Alliance on Systems Biology (SB Cancer, Submodule
V.7) and D.G. additionally by ViroQuant.

\bibliographystyle{abbrv}
\bibliography{prombib}

\end{document}